# A new complete algorithm for

# Irreducible Diophantine Pythagorean Triangles (IDPTs)

**JP Zwart   20220803**

Contents

# 1   Introduction

It is well known that a triangle with side lengths 3, 4 and 5 is right-angled. Euclid was the first to give a formula for generating other right-angled triangles with integer side lengths. Overviews of his and other approaches can be found in (consulted 26/7/2022):
https://en.wikipedia.org/wiki/Pythagorean_triple
https://en.wikipedia.org/wiki/Formulas_for_generating_Pythagorean_triples

In this text, I present a novel algorithm to generate all possible right-angled triangles with integer side lengths, in which the three side lengths have no common divisor. The algorithm is based on the difference in length between the hypothenuse and the largest of the two other sides. I also prove the completeness of this algorithm: it generates all possible such triangles and nothing but such triangles.

My term Diophantine Pythagorean Triangle (DPT) is synonymous with the term Pythagorean Triple used in the sources above, and my term Irreducible Diophantine Pythagorean Triangle (IDPT) is synonymous with Primitive Pythagorean Triple.

### a   Definition of IDPTs

A **Pythagorean Triangle** is a right-angled triangle.
Let $a$, $b$ and $c$ be the side lengths, with $a$ and $b$ the sides forming the right angle, and $c$ the hypothenuse. Then the sides satisfy Pythagoras' theorem: $a^2 + b^2 = c^2$. Obviously $c > a$ and $c > b$. It is also possible, without loss of generality, to choose $a \leq b$: otherwise just swap $a$ and $b$.

A **Diophantine Pythagorean Triangle (DPT)** is a right-angled triangle with integer side lengths.
For a DPT, $a$ cannot be equal to $b$, because then $c = \sqrt{2a^2} = \sqrt{2} \cdot a$, and so $c$ cannot be an integer.
So for a DPT: $a \in \mathbb{N}$, $b \in \mathbb{N}$, $c \in \mathbb{N}$, and $a < b < c$.
A DPT will be notated as a triplet $(a, b, c)$.
Familiar examples: (3, 4, 5) and (5, 12, 13).

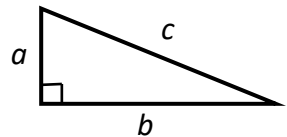

If $a$, $b$ and $c$ satisfy $a^2 + b^2 = c^2$, then obviously for any n $\in \mathbb{N}$: $(na)^2 + (nb)^2 = (nc)^2$ as well.
For example: $3^2 + 4^2 = 5^2$, so $(2\times3)^2 + (2\times4)^2 = (2\times5)^2$ $\Leftrightarrow$ $6^2 + 8^2 = 10^2$.
Therefore it is only interesting to consider numbers $a$, $b$ and $c$ that have no common factor in their prime factorizations. Since $6 = 2^1 \cdot 3^1$, $8 = 2^3$ and $10 = 2^1 \cdot 5^1$, there is a common factor 2 in their prime factorizations. In other words, if $a$, $b$ and $c$ have no common factor, then their Greatest Common Divisor is 1. This is notated here as **GCD($a$, $b$, $c$) =1**. If GCD($a$, $b$, $c$) =1 then $a$, $b$ and $c$ are **irreducible**.

An **Irreducible Diophantine Pythagorean Triangle (IDPT)** is a DPT in which GCD($a$, $b$, $c$) =1.
So the following holds:
**($a$, $b$, $c$) is an IDPT** $\Leftrightarrow$ $a \in \mathbb{N}$, $b \in \mathbb{N}$, $c \in \mathbb{N}$, $0 < a < b < c$, $a^2 + b^2 = c^2$ and **GCD($a$, $b$, $c$) =1**     [1]
Familiar examples of IDPTs are ($a$, $b$, $c$) = (3, 4, 5) and ($a$, $b$, $c$) = (5, 12, 13).

### b   Objective and outcome of this study

The objective of this monography is to find all possible IDPTs. The following chapters reflect the path I took to achieve this goal, and so might not offer the most efficient way to approach the problem.

The outcomes are: necessary and sufficient criteria for ($a$, $b$, $c$) to form an IDPT, and a complete algorithm for generating all IDPTs for $a \leq a$_max, with $a$_max an integer that can be chosen arbitrarily. There are infinitely many IDPTs, but any specific IDPT can be found using this algorithm.

## 2 Useful properties of square numbers.

Because $a^2 + b^2 = c^2$ can also be regarded in the form **$c^2 - b^2 = a^2$**, [2]
the properties of differences between squares might be useful.
Here are the first 16 positive integers with their squares,
and with $\Delta_n$ as the difference between the two adjacent square numbers $n^2$ and $(n+1)^2$.

| n | 1 | 2 | 3 | 4 | 5 | 6 | 7 | 8 | 9 | 10 | 11 | 12 | 13 | 14 | 15 | 16 |
|---|---|---|---|---|---|---|---|---|---|---|---|---|---|---|---|---|
| n² | 1 | 4 | 9 | 16 | 25 | 36 | 49 | 64 | 81 | 100 | 121 | 144 | 169 | 196 | 225 | 256 |
| Δn |  | 3 | 5 | 7 | 9 | 11 | 13 | 15 | 17 | 19 | 21 | 23 | 25 | 27 | 29 | 31 |

$\Delta_n$ can be calculated easily from **$\Delta_n = (n+1)^2 - n^2 = 2n+1$**. [3]

This example illustrates the following general property:
**The difference between two arbitrary squares is the sum of a row of consecutive odd numbers.**

This sum of consecutive numbers can be expressed in terms of $b$ and $c$ as follows:
 Let $b, c, d \in \mathbb{N}$, with **$c - b = d \Leftrightarrow b = c - d, d > 0$.** [4]
 Then $c^2 - b^2 = c^2 - (c-d)^2 = \{c^2-(c-1)^2\} + \{(c-1)^2-(c-2)^2\} + \ldots + \{(c-d+1)^2-(c-d)^2\} = \Delta_{c-1} + \Delta_{c-2} + \ldots + \Delta_{c-d}$.
 $c^2 - b^2 = \sum_{i=1}^{d}[(c-i+1)^2 - (c-i)^2] = \sum_{i=1}^{d} \Delta_{c-i} = \sum_{i=1}^{d}[2(c-i)+1]$
Likewise, alternatively:
 $c^2 - b^2 = (b+d)^2 - b^2 = \{(b+1)^2-b^2\} + \{(b+2)^2-(b+1)^2\} + \ldots + \{(b+d)^2-(b+d-1)^2\} = \Delta_b + \Delta_{b+1} + \ldots + \Delta_{b+d-1}$.
 $c^2 - b^2 = \sum_{i=1}^{d}[(b+i)^2 - (b+i-1)^2] = \sum_{i=1}^{d} \Delta_{b+i-1} = \sum_{i=1}^{d}[2((b-1)+i)+1]$

Thus proves the following property:
**Property 1** $\qquad c^2 - b^2 = \sum_{i=1}^{d}[2((b-1)+i)+1]$ [5a]
$\qquad\qquad\qquad c^2 - b^2 = \sum_{i=1}^{d}[2(c-i)+1]$ [5b]

Property 1 can be linked to DPTs as follows.
It uses $d$, the difference in length between the two largest sides, on the one hand ([4]: $d = c - b$),
and $c^2 - b^2$, the square of the length of the shortest side, on the other hand ([2]: $a^2 = c^2 - b^2$).
$\qquad a^2 = \sum_{i=1}^{d}[2((b-1)+i)+1]$ [6a]
$\qquad a^2 = \sum_{i=1}^{d}[2(c-i)+1]$ [6b]
[6a] and [6b] both state that $a^2$ is the sum of a row of length $d$ of consecutive odd numbers,
with [6a] linking $a$ to $b$ and $d$, and [6b] linking $a$ to $c$ and $d$.

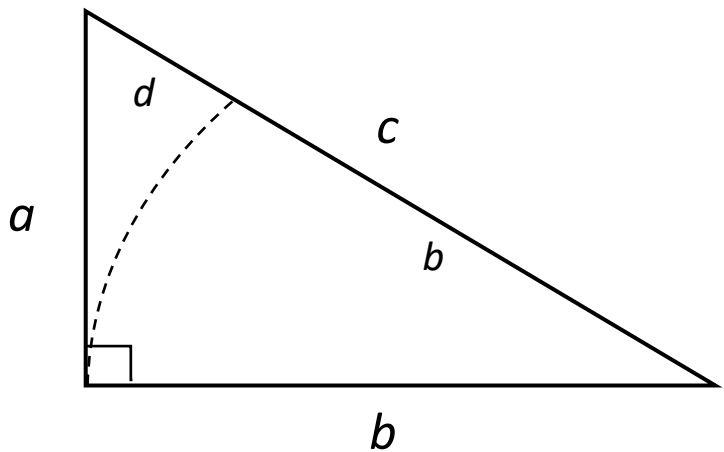

*Figure 1: DPT with lengths a, b, c and d*

# 3    First family of IDPTs

The simplest case in [6a] and [6b] is: **d = 1**. Then:

$$a^2 = 2b + 1 \Leftrightarrow b = \tfrac{1}{2}(a^2 - 1) \qquad [7a]$$
$$a^2 = 2c - 1 \Leftrightarrow c = \tfrac{1}{2}(a^2 + 1) \qquad [7b]$$

Because $b$ and $c$ must be integers greater than 0, $a^2$ must be an odd number greater than 1, and so $a$ can be any odd number greater than 1. Writing $a$ as $2n+1$, for $n \in \mathbb{N}$, $n > 0$, this becomes:

$$a = 2n + 1 \qquad [8a]$$
$$b = \tfrac{1}{2}(a^2 - 1) = \tfrac{1}{2}((2n+1)^2 - 1) = 2n^2 + 2n \qquad [8b]$$
$$c = \tfrac{1}{2}(a^2 + 1) = \tfrac{1}{2}((2n+1)^2 + 1) = 2n^2 + 2n + 1 \qquad [8c]$$

This proves the following statement:

**Lemma 1**

**For any number $n \in \mathbb{N}$, $n > 0$, there exists a DPT ($a$, $b$, $c$) = ($2n+1$, $2n^2+2n$, $2n^2+2n+1$).** [9]

Examples:  ( 3,  4,  5) with    9 +   16 =   25
          ( 5, 12, 13) with   25 +  144 =  169
          ( 7, 24, 25) with   49 +  576 =  625
                …
          (13, 84, 85) with  169 + 7056 = 7225
                …

Do any of these triplets ($a$, $b$, $c$) have a common divisor, or are they irreducible?
Observe, that $c - b = 2n^2+2n+1 - (2n^2+2n) = 1$.

**Lemma 2: If $c - b = 1$, then $c$ and $b$ have no common divisor.** [10]

Proof:  Every integer greater than 1 has a unique factorization in primes.
        So $b = p_1 \cdot p_2 \cdot \ldots \cdot p_n$, for some number n and with primes $p_i \leq p_{i+1}$ for $1 \leq i \leq n-1$.
        Then $c = b + 1 = p_1 \cdot p_2 \cdot \ldots \cdot p_n + 1$.
        So if $c$ is divided by any prime $p_i$ from the factorization of $b$, the remainder is always 1.
        Therefore $c$ is not divisible by any of the divisors of $b$, so they have no divisor in common.

Lemma 2 shows that $b$ and $c$ have no common divisor, so all DPTs from [9] are irreducible.

**Proposition 1**

**Let $d = c - b$. If $d = 1$, then for any number $n \in \mathbb{N}$, $n > 0$,**
**there exists an IDPT ($a$, $b$, $c$) = ($2n+1$, $2n^2+2n$, $2n^2+2n+1$).** [11]

Proof: follows directly from Lemma 1 in [9] and lemma 2 in [10].

**Conclusions**

Proposition 1 in [11] defines the **first family of IDPTs**. There are countably infinite different IDPTs in this family. No other IDPTs with $c - b = 1$ exist: all possible values of $a$ have been considered. This family was found by considering [6a] and [6b] with $d = 1$. Perhaps other families can be found with $d = 2, 3, \ldots$. The following chapters will explore this.

## 4   Expressing *a*, *b* and *c* in terms of *d* and *n*, for any *d* > 0

In an arbitrary triangle, there are three degrees of freedom in choosing the side lengths *a*, *b* and *c*. Two of these lengths can be chosen completely freely as $l_1$ and $l_2$, the third can then be any number between $(l_1+l_2)$ and $|l_1-l_2|$. But if the form of the triangle is restricted to being right-angled, then *c* is fixed once *a* and *b* are chosen, which leaves only two degrees of freedom. Because the obvious restriction $c = \sqrt{a^2 + b^2}$ is inconvenient for finding integer solutions, expressions for *a*, *b* and *c* in terms of two integer variables are derived below. One of these is the difference *d* between *c* and *b*. For *d* = 1, chapter 3 found a family of IDPTs. This chapter generalizes this result for $d \geq 1$.

*a*, *b* and *c* must satisfy:  $\qquad\qquad\qquad\qquad$ **$a^2 + b^2 = c^2$, with $0 < a < b < c$** $\qquad\qquad$ [1]
*d* is defined as $d := c - b$: $\qquad\qquad\qquad\qquad$ **$c = b + d$, with d > 0** $\qquad\qquad\qquad\qquad$ [5]
$a^2 + b^2 = c^2 \Leftrightarrow a^2 = c^2 - b^2 = (b + d)^2 - b^2 = 2bd + d^2 \Leftrightarrow$ **$b = (a^2 - d^2)/2d$** $\qquad$ [14]
$c = b + d$ $\qquad\qquad\qquad\qquad\qquad\qquad\Leftrightarrow$ **$c = (a^2 - d^2)/2d + d = (a^2 + d^2)/2d$** $\qquad$ [15]
So *b* and *c* must be integers that satisfy [1], [14] and [15], with *d* the distance between *b* and *c*.
It follows from [1] and [14] that $a^2$ must be greater than $d^2$, and so also ***a > d***. $\qquad\qquad$ [16]
A DPT is given by **$(a, b, c) = (a, (a^2 - d^2)/2d, (a^2 - d^2)/2d + d)$, with a > d** $\qquad$ [17]
An IDPT is given by [17] with the additional requirement that **GCD(*a*, *b*, *c*) = 1** $\qquad$ [18]

Because of [14], since *b* must be an integer, $(2d \mid (a^2 - d^2))$.
A simpler expression to find numbers that satisfy this divisibility can be obtained as follows.
From [17]: $b = (a^2 - d^2)/2d = (a - d)(a + d)/2d$, with $a > d$.
Let *m* > 0 be defined as: $\qquad\qquad\qquad\qquad$ **$m := a - d$** $\qquad\qquad\qquad\qquad\qquad\qquad$ [19]
Using [19] for *a* yields: $\qquad\qquad\qquad\qquad$ **$a = m + d$, with $m > 0$, $d > 0$** $\qquad\qquad$ [20]
Substituting [19] into [14]: $b = m(m + 2d)/2d$ $\quad\Leftrightarrow$ **$b = m^2/2d + m$, with $m > 0$, $d > 0$** $\quad$ [21]
Using [5], a DPT is given by **$(a, b, c) = (m + d, m^2/2d + m, m^2/2d + m + d)$, with $m > 0$, $d > 0$** $\quad$ [22]
See figure 2, which shows *d* and *m* in a DPT. The meaning of *n* is explained below figure 2.

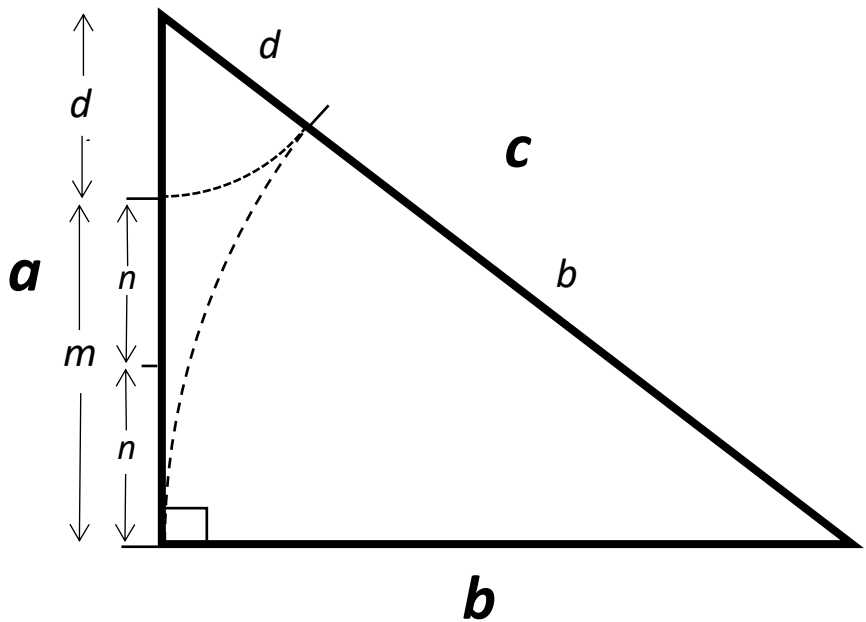

*Figure 2: DPT with d = c-b, m = a-d and m = 2n*

A further restriction on the minimum value of *m* can be obtained from [1]: *a* < *b* and [20]-[21]:
*m* + *d* < $m^2$/2*d* + *m* ⇔ 2$d^2$ < $m^2$ ⇔ **$m^2$ > 2$d^2$**  [23]
So [23] sets a higher minimum for *m* than [19]-[22]. Example: *d* = 5 ⇒ *m* ≥ 8.

Because of [21], since *b* must be an integer, (2*d* | $m^2$).
Therefore [21] and [22] can be written more conveniently as follows.
(2*d* | $m^2$) ⇒ ∃ g ∈ ℕ: 2*dg* = $m^2$ ⇒ (2 | $m^2$) ⇔ (2 | *m*).
In other words: $m^2$ is an even number, and so *m* is even as well, otherwise $m^2$ would be odd.

(2 | *m*) ⇒ ∃ *n* ∈ ℕ: **2*n* = *m***. Let *n* be defined as    ***n* := *m* / 2 ⇔ 2*n* = *m*.**  [24]
Substituting [24] into [20]-[23] yields    ***a* = 2*n* + *d*, with *d* > 0, 2$n^2$ > $d^2$**  [25]
***b* = 2$n^2$/*d* + 2*n*, with *d* > 0, 2$n^2$ > $d^2$**  [26]
A DPT is given by **(*a*, *b*, *c*) = (2*n* + *d*, 2$n^2$/*d* + 2*n*, 2$n^2$/*d* + 2*n* + *d*), with *d* > 0, 2$n^2$ > $d^2$, (*d* | 2$n^2$)**  [27]
An IDPT is given by [27] with the additional requirement that **GCD(*a*, *b*, *c*) = 1**  [28]
See figure 2 above, which shows *d* and *n* in a DPT.

It follows from [25], since *n* > 0, that *a* ≥ *d* + 2, though a higher minimum can follow from 2$n^2$ > $d^2$.
The following inequalities must all hold: **0 < *a* < *b* < *c*, *a* ≥ *d* + 2, 2$n^2$ > $d^2$**  [29]
The last part of [29] sets a minimum for *n* for each given *d*. Example: *d* = 5 ⇒ *n* ≥ 4.
Equations [25] – [29] will be used below to explore possible values of (*a*, *b*, *c*) for each given *d*.

The approach in the chapters that follow is to establish for each value of *d* whether IDPTs exist, and if so, what they are.

# 5    Integer term 2n²/d: cases *d* = 1 and *d* = 2

[27] states: (*a*, *b*, *c*) = ($2n + d$, $2n^2/d + 2n$, $2n^2/d + 2n + d$), with $d > 0$, $2n^2 > d^2$, ($d \mid 2n^2$).
Since *b* and *c* must be integers, *d* must be a divisor of $2n^2$: ($d \mid 2n^2$).
There are therefore two special cases: *d* = 1 and *d* = 2. In both these cases ($d \mid 2n^2$) for all $n \in \mathbb{N}$.
These cases are discussed in this chapter.
For *d* > 2, it is more difficult to establish for which $n \in \mathbb{N}$ ($d \mid 2n^2$) holds. These cases are discussed in chapter 6.

### a    Case *d* = 1
From [27] with *d* =1: (*a*, *b*, *c*) = ($2n + 1$, $2n^2 + 2n$, $2n^2 + 2n + 1$), $2n^2 > 1 \Leftrightarrow n \geq 1$.
This yields the first family of IDPTs, already discussed in chapter 3.

### b    Case *d* = 2
From [27] with *d* =2: (*a*, *b*, *c*) = ($2n + 2$, $n^2 + 2n$, $n^2 + 2n + 2$), $2n^2 > 4 \Leftrightarrow n \geq 2$.
However, if ($2 \mid n$) then *a*, *b* and *c* have a common factor 2,
because ($2 \mid n$) $\Leftrightarrow$ ($2 \mid n^2$) and so every term in ($2n + 2$, $n^2 + 2n$, $n^2 + 2n + 2$) is divisible by 2.
Therefore *n* must de an odd number.

**Lemma 3**
**For any odd number *n* > 1, there exists a DPT (*a*, *b*, *c*) = ($2n + 2$, $n^2 + 2n$, $n^2 + 2n + 2$).**    [30]
Examples: ( 8, 15,  17) with  64 +  225 =   289
         (12, 35,  37) with 144 + 1225 =  1369
         (16, 63,  65) with 256 + 3969 =  4225
         (20, 99, 101) with 400 + 9801 = 10201
              …
Proof:  it is easily verified algebraically that $a^2 + b^2 = c^2$ with (*a*, *b*, *c*) given by [27],
        for any value of *n* > 0 and *d* > 0, so in particular this holds for *d* = 2, and for *n* odd, *n* > 1.

**Lemma 4: If *c* is odd and *c* – *b* = 2, then *c* and *b* have no common factor.**    [31]
Proof:  Every integer greater than 1 has a unique factorization in primes.
        So $b = p_1 \cdot p_2 \cdot \ldots \cdot p_n$, for some number n and with primes $p_i \leq p_{i+1}$ for $1 \leq i \leq$ n-1.
        Then $c = b + 2 = p_1 \cdot p_2 \cdot \ldots \cdot p_n + 2$.
        Because *c* is odd, *b* = *c* – 2 is odd as well, and 2 is not a divisor of *b*.
        So if *c* is divided by any prime $p_i > 2$ from the factorization of *b*, the remainder is always 2.
        Therefore *c* is not divisible by any of the divisors of *b*, so they have no divisor in common.
Lemma 4 shows that *b* and *c* have no common divisor, so all DPTs from [30] are irreducible.

**Proposition 2**
**Let *d* = *c* – *b*. If *d* = 2, then for any odd number *n* greater than 1,**
**there exists an IDPT (*a*, *b*, *c*) = ($2n + 2$, $n^2 + 2n$, $n^2 + 2n + 2$).**    [32]
Proof: follows directly from Lemma 3 in [30] and lemma 4 in [31].

**Conclusions:**
Proposition 2 in [32] defines the **second family of IDPTs**. There are countably infinite different
IDPTs in this family. No other IDPTs with *d* = 2 exist: all possible values of *n* have been considered.

# 6   Integer term $2n^2/d$: cases $d > 1$

To find more values of *d* for which IDPTs exist, the prime factorization of *d* is to be considered. The notation used in this text for prime factorizations is presented in appendix A, part d.

## a   Conditions ($d \mid 2n^2$) and GCD(*a*, *b*, *c*) = 1

A DPT is given by **(*a*, *b*, *c*) = ($2n + d$, $2n^2/d + 2n$, $2n^2/d + 2n + d$)**, with $d > 0$, $2n^2 > d^2$, ($d \mid 2n^2$)   [27]

For *b* and *c* to be integers, a necessary and sufficient condition for *d* and *n* is: **($d \mid 2n^2$)**.   [33]

But condition [33] is not enough for (*a*, *b*, *c*) to form an irreducible DPT:

for (*a*, *b*, *c*) given by [27] to be an IDPT, also **GCD(*a*, *b*, *c*) = 1** must hold.   [34]

It turns out, that a useful necessary condition for [34] to hold can be derived. It will be stated and proved in section 6b below. This condition will produce all the possible values of *d* for which IDPTs exist. There are not many such values for *d*: below 100, *d* can only have the values 1, 2, 8, 9, 18, 25, 32, 49, 50, 72, 81 and 98, for example.

Still, the condition of section 6b is not sufficient: for each possible value of *d* there might be some values of *n* that indeed yield IDPTs, and other values of *n* that do not. This was the case in the second family of IDPTs described in section 5b, for example, in which only odd values of n yielded IDPTs. So for each possible value of *d*, further restrictions must be imposed on *n*, the remaining degree of freedom, to comply with [34]. The discussion on possible values of *n* depending on the value of *d* starts in chapter 7. The final aim is to arrive at necessary and sufficient conditions for [34] to hold.

## b   Necessary condition on $d > 1$

**A necessary condition on $d > 1$ for [34] is: in PF(*d*), $e_2(d)$ is odd, and $\forall_{p>2}$: $e_p(d)$ is even.**   [35]

In this text, this condition will also be called **condition NC1.**   [36]

To prove this, common divisors of *a*, *b* and *c*, and of all the terms in [27] must be considered.

This is done using the prime factorizations of all these numbers, presented in appendix A, part d.

This can be stated as a theorem:

**Theorem 1**

**(*a*, *b*, *c*) is an IDPT with d > 1** $\Rightarrow$ ( $\neg(2 \mid e_2(d)) \wedge \forall_{p>2}$: ($2 \mid e_p(d)$))   [37]

The proof is presented in the three subsections below.

### i   Necessary condition on $e_2(d)$

**Lemma 5**

**If (*a*, *b*, *c*) is a DPT from [27], and $e_2(d)$ is an even number > 0, then (*a*, *b*, *c*) is reducible.**   [38]

Proof:   Suppose **$e_2(d) = 2k$, $k \in \mathbb{N}$, $k > 0$**. Then **($2 \mid d$)**.   [39]

Then all the terms in [27]: (*a*, *b*, *c*) = ($2n + d$, $2n^2/d + 2n$, $2n^2/d + 2n + d$), are divisible by 2, provided that $2n^2/d$ is also divisible by 2.

Note: in general it is not true that ($2 \mid 2p/q$), but only if p/q is itself an integer.

(example: $\neg(2 \mid (2 \cdot 6^2)/8)$). It must therefore be explicitly shown that ($2 \mid 2n^2/d$).

In an IDPT, ($d \mid 2n^2$), so $2n^2/d \in \mathbb{N}$.   [40]

From [40], PF($2n^2/d$) and properties [PP1] - [PP2] follows:

$e_2(2n^2/d) = 2e_2(n) + 1 - e_2(d) \geq 0$   [41]

However, can $e_2(2n^2/d) = 0$?

No, because then, from [39] and [41]: $2e_2(n) + 1 - e_2(d) = 0 \Leftrightarrow e_2(n) = \frac{1}{2}(e_2(d) - 1) = k - \frac{1}{2}$.

This contradicts property P1: $e_2(n)$ must be an integer.

So $e_2(2n^2/d) = 2e_2(n) + 1 - e_2(d) > 0 \Leftrightarrow e_2(2n^2/d) > 0 \Rightarrow$ **($2 \mid 2n^2/d$)**.   [42]

It follows from [39] and [42] that all terms in [27] are even, and so *a*, *b* and *c* have a common divisor 2. So: $(2 \mid e_2(d)) \Rightarrow$ GCD(a, b, c) $\geq 2 > 1$, contradicting condition [34]: GCD(a, b, c)=1.
So $e_2(d)$ cannot be an even number > 0:

**$(2 \mid e_2(d) \Rightarrow \neg$IDPT(a, b, c) with *d* > 1] $\Leftrightarrow$ IDPT(a, b, c) with *d* > 1 $\Rightarrow \neg(2 \mid e_2(d)$**  [43]

A necessary condition on *d* for (a, b, c) to be an IDPT with d > 1 is therefore: **$\neg(2 \mid e_2(d))$**  [44]

      **ii**      **Necessary condition on $e_p(d)$, for all p > 2.**

**Lemma 6**

**If (*a*, *b*, *c*) is a DPT from [27], and $\exists$ p>2: $e_p(d)$ is an odd number > 0, then (*a*, *b*, *c*) is reducible.** [45]

Proof:   Suppose $\exists$ **p > 2: $e_p(d)$ = 2k+1, k $\in \mathbb{N}$, k $\geq$ 0**. Then **(p | *d*)**.  [46]

If (p | *d*), all the terms in (*a*, *b*, *c*) = $(2n + d, 2n^2/d + 2n, 2n^2/d + 2n + d)$ are divisible by p, provided that $2n^2/d$ and 2n are also divisible by p.

In an IDPT, $(d \mid 2n^2)$, so $2n^2/d \in \mathbb{N}$.  [47]

From [47], PF($2n^2/d$) and properties [P1] - [P3] follows: $e_p(2n^2/d) = 2e_p(n) - e_p(d) \geq 0$  [48]

However, can $e_p(2n^2/d) = 0$?

No, because then, from [46] and [48]: $2e_p(n) - e_p(d) = 0 \Leftrightarrow e_p(n) = \tfrac{1}{2}e_p(d) = k + \tfrac{1}{2}$.

This contradicts property P1: $e_p(n)$ must be an integer.

So $e_p(2n^2/d) = 2e_p(n) - e_p(d) > 0 \Leftrightarrow e_p(2n^2/d) > 0 \Rightarrow$ **(p | $2n^2/d$)**.  [49]

Is p also a divisor of 2n? Yes, because:

From [46]: $\exists$ g $\in \mathbb{N}$, g > 0: $d/p = g \Leftrightarrow$ ***d* = gp**  [50]

From [49]: $\exists$ h $\in \mathbb{N}$, h > 0: $(2n^2/d)/p = h \Leftrightarrow$ **$2n^2/d$ = hp**  [51]

From [50] and [51]: $2n^2/d = 2n^2/gp = hp \Leftrightarrow 4n^2/p^2 = 2gh \Leftrightarrow$ **$(2n/p)^2$ = 2gh**  [52]

So $(2n/p)^2 \in \mathbb{N}$. 2n/p itself is rational because n $\in \mathbb{N}$ and p $\in \mathbb{N}$. But the square of a rational number can only be an integer if it is an integer itself. So 2n/p $\in \mathbb{N} \Rightarrow$ **(p | 2n)**  [53]

It follows from [46], [49] and [53] that p is a divisor of all the terms in [27], and so *a*, *b* and *c* have a common divisor p. So: $\neg(2 \mid e_p(d)) \Rightarrow$ GCD(a, b, c) $\geq$ p > 1, contradicting [34].

So $e_p(d)$ cannot be an odd number > 0:

**$\neg(2 \mid e_p(d)) \Rightarrow \neg$IDPT(a, b, c) with *d* > 1 $\Leftrightarrow$ IDPT(a, b, c) with *d* > 1 $\Rightarrow (2 \mid e_2(d))$**  [54]

A necessary condition on *d* for (a, b, c) to be an IDPT with d > 1 is thus: **$\forall_{p>2}$: (2 | $e_p(d)$)]**  [56]

      **iii**      **Combination of conditions on $e_p(d)$: necessary condition on *d***

Lemma 5 proves that condition [44] is a necessary condition on $e_2(d)$, and lemma 6 proves that condition [56] is a necessary condition on all the other prime divisors in PF(*d*). Together, [44] and [56] cover all the prime divisors in PF(*d*). This proves theorem 1 in [36]:

**Theorem 1:**

**(a, b, c) is an IDPT with *d* > 1 $\Rightarrow \neg(2 \mid e_2(d))$} $\wedge \forall_{p>2}$: (2 | $e_p(d)$)**  [37]

or in shorthand:

**(a, b, c) is an IDPT with *d* > 1 $\Rightarrow$ NC1**  [38]

This means: values of *d* that do not meet NC1 will certainly not yield an IDPT. NC1 can be extended to exclude cases that satisfy NC1 yet do not form an IDPT; in other words: NC1 can be extended so it will also be a sufficient condition for (a, b, c) to be an IDPT. It will be shown in chapters to follow, that for each value if *d* that meets NC1 in [38], restrictions on the values of *n* can be found such that (a, b, c) is indeed an IDPT. These restrictions on *n*, together with NC1, then form also a sufficient condition for (a, b, c) to be an IDPT.

# 7 Cases $d = 2 \cdot 2^{2k}$, k ≥ 0: generalization of the second family of IDPTs

Condition NC1 on $d > 1$ states: $\neg(2 \mid e_2(d)) \land \forall_{p>2}: (2 \mid e_p(d))$ [57]

This chapter considers only the values of $d > 1$ that are a power of 2: $e_2(d) > 0$. Other values of $d$ are discussed in subsequent chapters. A general formula for all possible IDPTs with $d$ being a power of 2 is derived in section 7a. This formula is a generalization of the second family of IDPTs, which is described in section 5a above. The proof of the formula in 7a needs some properties of divisibility, which are presented in appendix A, part e. Section 7b summarizes the results of section 7a in a theorem. Section 7c finally presents examples of IDPTs for several values of $d$ and $n$.

## a  Values of *n* for $d = 2 \cdot 2^{2k}$, k ≥ 0

This section deals with all values of $d$ that are odd powers of 2: $d = 2^{(2k+1)} = 2 \cdot 2^{2k}$. The second family of IDPTs, discussed in section 5b above, falls within this set of $d$-values, with k = 0. For each value of $d$, values of $n$ are given such that all resulting triangles are IDPTs. The entire set of $d$-values together with their corresponding $n$-values makes up the extended second family of IDPTs.

### i   DPTs with $d = 2 \cdot 2^{2k}$, k ≥ 0

A DPT is given by **(a, b, c) = (2n + d, 2n²/d + 2n, 2n²/d + 2n + d)**, with $d > 0$, $2n^2 > d^2$, $(d \mid 2n^2)$ [27]

With $d = 2 \cdot 2^{2k}$, the requirement $(d \mid 2n^2)$, needed for $b$ and $c$ to be integers, leads to:
$2n^2/(2 \cdot 2^{2k}) \in \mathbb{N} \Leftrightarrow n^2/2^{2k} \in \mathbb{N} \Leftrightarrow (n/2^k)^2 \in \mathbb{N} \Leftrightarrow n/2^k \in \mathbb{N} \Leftrightarrow n/2^k = r \in \mathbb{N} \Leftrightarrow$
**$n = 2^k \cdot r, r \in \mathbb{N}, r \geq 1$.** [58]

With $d = 2 \cdot 2^{2k}$ and $n = 2^k \cdot r$, the requirement $2n^2 > d^2$, needed for $a$ to be smaller than $b$, leads to:
$2 \cdot (2^k \cdot r)^2 > d^2 \Leftrightarrow 2 \cdot 2^{2k} \cdot r^2 > d^2 \Leftrightarrow d \cdot r^2 > d^2 \Leftrightarrow r^2 > d \Leftrightarrow r > \sqrt{d} \Leftrightarrow$
**$r \geq \text{ceiling}(\sqrt{d}) = r_{min,2,k}$.** [59]

The prime divisor 2 is added to the index of $r_{min,2,k}$ because other primes will be considered later as well, and the general notation $r_{min, p, k}$ will be used for any prime divisor p.

Table T1 lists values of r and $n$ for a few values of k. It can be easily extended.

| k | $d = 2 \cdot 2^{2k}$ | $r_{min,2,k}$ | r | $n = 2^k \cdot r$ |
|---|---|---|---|---|
| 0 | 2 | 2 | 2, 3, 4, 5, … | 2, 3, 4, 5, … |
| 1 | 8 | 3 | 3, 4, 5, 6, … | 6, 8, 10, 12, … |
| 2 | 32 | 6 | 6, 7, 8, 9, … | 24, 28, 32, 36, … |
| 3 | 128 | 12 | 12, 13, 14, 15, … | 96, 104, 112, 120, … |
| 4 | 512 | 23 | 23, 24, 25, 26, … | 368, 384, 400, 416, … |
| … | … | … | … | … |

*Table T1: values of n to form a DPT for given $d = 2 \cdot 2^{2k}$*

**$n = 2^k \cdot r$, $d = 2 \cdot 2^{2k}$, $2n^2/d = 2(2^k \cdot r)^2/2^{(2k+1)} = r^2$** and **$r = r_{min,2,k}, r_{min,2,k} + 1, r_{min,2,k} + 2, …$** [60r]

With [60r], the DPTs ($a, b, c$) are given by:
    $a = 2n + d$         $= \mathbf{2 \cdot 2^k \cdot r + 2 \cdot 2^{2k}}$ [60a]
    $b = 2n + r^2$        $= \mathbf{2 \cdot 2^k \cdot r + r^2}$ [60b]
    $c = = a + r^2 = b + d$  $= \mathbf{2 \cdot 2^k \cdot r + 2 \cdot 2^{2k} + r^2}$ [60c]

### ii      GCD($a$, $b$, $c$) > 1 for DPTs with $d = 2 \cdot 2^{2k}$ and (2 | r)

For a DPT ($a$, $b$, $c$) given by [60a] – [60c], and r given by [60r], **GCD($a$,$b$,$c$) > 1 if r is an even number**.
    (2 | r) $\Rightarrow$ (2 | $a$), because all the terms in [60a] are even.
    (2 | r) $\Rightarrow$ (2 | $b$), because all the terms in [60b] are even.
    (2 | r) $\Rightarrow$ (2 | $c$), because all the terms in [60c] are even.
So if (2 | r), then GCD($a$, $b$, $c$) $\geq$ 2 and DPT ($a$, $b$, $c$) is reducible.

### iii      IDPTs with $d = 2 \cdot 2^{2k}$, $k \geq 0$

This section removes those DPTs from section 7ai, which result from an even value of r. It will be proved in section 7aiv, that the remaining DPTs are indeed IDPTs.

As is shown in section 7aii: (2 | r) $\Rightarrow$ GCD($a$, $b$, $c$) > 1.
Therefore a necessary condition for DPT ($a$, $b$, $c$) to be irreducible, is: ¬(2 | r)           [61]

Incorporating [61] into table T1 yields table T2, which lists values of r and $n$ for a few values of k. It can be easily extended.

| k | $d = 2 \cdot 2^{2k}$ | $r_{min,2,k}$ | r | $n = 2^k \cdot r$ |
|---|---|---|---|---|
| 0 | 2 | 2 | 3, 5, 7, 9, … | 3, 5, 7, 9, … |
| 1 | 8 | 3 | 3, 5, 7, 9, … | 6, 10, 14, 18, … |
| 2 | 32 | 6 | 7, 9, 11, 13, … | 28, 36, 44, 52, … |
| 3 | 128 | 12 | 13, 15, 17, 19, … | 104, 120, 136, 152, … |
| 4 | 512 | 23 | 23, 25, 27, 29, … | 368, 400, 432, 464, … |
| … | … | … | … | … |

*Table T2: values of n to form an IDPT for given $d = 2 \cdot 2^{2k}$*

$d = 2^{(2k+1)}$,   $n = 2^k \cdot r$,   $2n^2/d = r^2$   and   $r \geq r_{min,2,k}$,   ¬(2 | r)         [62n]
With [62n], the IDPTs ($a$, $b$, $c$) are given by:
    $a = 2n + d$         = **$2 \cdot 2^k \cdot r + 2 \cdot 2^{2k}$**         [62a]
    $b = 2n + r^2$       = **$2 \cdot 2^k \cdot r + r^2$**         [62b]
    $c = = a + r^2 = b + d$   = **$2 \cdot 2^k \cdot r + 2 \cdot 2^{2k} + r^2$**     [62c]

A few examples of IDPTs from T2 and [62]:
    k = 0, $d$ = 2, r = 7, $n$ = 7:         ($a$, $b$, $c$) = (16, 63, 65)
    k = 1, $d$ = 8, r = 5, $n$ = 10:       ($a$, $b$, $c$) = (28, 45, 53)
    k = 2, $d$ = 32, r = 13, n = 52     ($a$, $b$, $c$) = (136, 273, 305)

### iv    Proof that the DPTs from [62] are irreducible

It will be shown, that no prime can be a common divisor of *a*, *b*, and *c* given by [62].
Then GCD(*a*, *b*, *c*) = 1, and (*a*, *b*, *c*) is an IDPT.
Each step below uses the obvious property: $\neg(g \mid x) \Rightarrow g \notin \{CD(x, y)\}$ for any numbers x and y.
See [ND2] in appendix A for the notation used.

**Step 1: $2 \notin \{CD(a, b, c)\}$**

    [62b]: $b = 2 \cdot 2^k \cdot r + r^2$.
    p = 2 $\Rightarrow$ p is a divisor of the first term in [62b]: $2 \cdot 2^k \cdot r$
    p = 2 $\Rightarrow$ p is not a divisor of the second term in [62b]: $r^2$,
    because $\neg(2 \mid r)$, and using [PD1b] from appendix A: $\neg(2 \mid r) \Rightarrow \neg(2 \mid r^2)$.
    Using [PD4a] from appendix A: $(2 \mid 2 \cdot 2^k \cdot r) \wedge \neg(2 \mid r^2) \Rightarrow \neg(2 \mid b)$
    So $2 \notin \{CD(a, b, c)\}$

**Step 2: Suppose p is a prime, p > 2 and (p | r). Then $p \notin \{CD(a, b, c)\}$**

    [62a]: $a = 2 \cdot 2^k \cdot r + 2 \cdot 2^{2k}$
    $(p \mid r) \Rightarrow$ p is a divisor of the first term in [62a]: $2 \cdot 2^k \cdot r$.
    p > 2 $\Rightarrow$ p is not a divisor of the second term in [62a]: $2 \cdot 2^{2k}$.
    Using [PD4a] from appendix A: $(p \mid 2 \cdot 2^k \cdot r) \wedge \neg(p \mid 2 \cdot 2^{2k}) \Rightarrow \neg(p \mid a)$
    So $p \notin \{CD(a, b, c)\}$

**Step 3: Suppose p is a prime, p > 2 and $\neg(p \mid r)$. Then $p \notin \{CD(a, b, c)\}$**

    **Step 3a: Suppose p is a prime, p > 2, $\neg(p \mid r)$ and also $\neg(p \mid a)$.**
        Then $p \notin \{CD(a, b, c)\}$.
    **Step 3b: Suppose p is a prime, p > 2, $\neg(p \mid r)$ and also $(p \mid a)$.**
        [62c]: $c = a + r^2$.
        $(p \mid a) \Rightarrow$ p is a divisor of the first term in [62c]: *a*.
        $\neg(p \mid r)$: Using [PD1b] from appendix A: $\neg(p \mid r) \Rightarrow \neg(p \mid r^2)$.
        So p is not a divisor of the second term in [62c]: $r^2$.
        Using [PD4a] from appendix A: $(p \mid a) \wedge \neg(p \mid r^2) \Rightarrow \neg(p \mid c)$
        So $p \notin \{CD(a, b, c)\}$

Steps 1, 2, 3a and 3b cover all possible prime divisors,
therefore $\forall_p: p \notin \{CD(a, b, c)\} \Leftrightarrow GCD(a, b, c) = 1$
This completes the proof that all DPTs given by [62] are irreducible.

## b    Extended second family of IDPTs

The extended second family of IDLTs is given by theorem 2 below.
It concerns all IDPTs in which *d* is an odd power of 2.

**Theorem 2**
Let *d*, r and *n* satisfy criterion C2, consisting of C2d, C2r and C2n:

$\quad$ [C2d]: $d = 2 \cdot 2^{2k}$, $k \in \mathbb{N}$, $k \geq 0$. $\hfill$ [C2d]

$\quad$ [C2r]: $r_{min,2,k} = \text{ceiling}(\sqrt{d})$. $\quad r \in \mathbb{N}$: $(r \geq r_{min,2,k}) \wedge \neg(2 \mid r)$. $\hfill$ [C2r]

$\quad$ [C2n]: $n = 2^k \cdot r$. $\hfill$ [C2n]

Let triangles (*a*, *b*, *c*) be given by:

$\quad a = 2n + d \qquad\qquad = 2 \cdot 2^k \cdot r + 2 \cdot 2^{2k}$ $\hfill$ [62a]

$\quad b = 2n + 2n^2/d \qquad = 2 \cdot 2^k \cdot r + r^2$ $\hfill$ [62b]

$\quad c = a + 2n^2/d = b + d \quad = 2 \cdot 2^k \cdot r + r^2 + 2 \cdot 2^{2k}$ $\hfill$ [62c]

Then all triangles (*a*, *b*, *c*) are IDPTs, and no other IDPTs with *d* satisfying C2d exist.
Proof: section 7a proves this theorem.

All these IDPTs form the **extended second family of IDPTs.**
The second family of IDPTs, discussed in section 5b, is the case in theorem 2 with k = 0.

## c    Examples of IDPTs in the extended second family

**k = 0, *d* =2**
**s$_{min,k}$ = 1, *n* = 3, 5, 7, 9, …**
IDPTs: (8, 15, 17), (12, 35, 37), (16, 63, 65), (20, 99, 101), …

**k = 1, *d* =8**
**s$_{min,k}$ = 1, *n* = 6, 10, 14, 18, …**
IDPTs: (20, 21, 29), (28, 45, 53), (36, 77, 85), (44, 117, 125), …

**k = 2, *d* =32**
**s$_{min,k}$ = 3, *n* = 28, 36, 44, 52, …**
IDPTs: (88, 105, 137), (104, 153, 185), (120, 209, 241), (136, 273, 305), …

**k = 3, *d* =128**
**s$_{min,k}$ = 6, *n* = 104, 120, 136, 152, …**
IDPTs: (336, 377, 505), (368, 465, 593), (400, 561, 689), (432, 665, 793), …

**k = 4, *d* =512**
**s$_{min,k}$ = 11, *n* = 368, 400, 432, 464, …**
IDPTs: (1248, 1265, 1777), (1312, 1425, 1937), (1376, 1593, 2105), (1440 1769 2281), …

**k = 5, *d* =2048**
**s$_{min,k}$ = 23, *n* = 1504, 1568, 1632, 1696, …**
IDPTs: (5056, 5217, 7265), (5184, 5537, 7585), (5312, 5865, 7913), (5440, 6201, 8249), …

….

# 8 Cases *d* = p$^{2k}$, k ≥ 0: generalization of the first family of IDPTs

Condition NC1 on *d* > 1 states: ¬(2 | e$_2$(*d*))} ∧ ∀$_{p>2}$: (2 | e$_p$(*d*))]}  [37]

Although all exponents in PF(*d*) must be greater than 0, the case k = 0 will also be considered here, because this allows the case *d* = 1 to be included. In chapter 9, strict conditions for this will be given. This chapter considers only the values of d > 1 that are an even power of a single prime number p, p>2: e$_p$(*d*) = 2k, k ≥ 0 and ∀$_{q \neq p}$: (e$_q$(*d*) = 0). Other values of d are discussed in subsequent chapters. A general formula for all possible IDPTs with *d* being an even power of p is derived in section 8a. This formula is a generalization of the first family of IDPTs, which is described in chapter 3 above. Section 8b summarizes the results of section 8a in a theorem. Section 8c finally presents examples of IDPTs for several values of *d* and *n*.

## a     Values of *n* for *d* = p$^{2k}$, k ≥ 0

This section deals with all values of *d* that are even powers of a prime number p, p > 2. The first family of IDPTs, discussed in chapter 3 above, falls within this set of *d*-values, with k = 0, for any p. For each value of *d*, values of *n* are given such that all resulting triangles are IDPTs. The entire set of *d*-values together with their corresponding *n*-values makes up the extended first family of IDPTs.

### i     DPTs with *d* = p$^{2k}$, k ≥ 0

A DPT is given by **(a, b, c) = (2n + d, 2n$^2$/d + 2n, 2n$^2$/d + 2n + d)**, with *d* > 0, 2n$^2$ > d$^2$, (d | 2n$^2$)   [27]

With *d* = p$^{2k}$, the requirement **(d | 2n$^2$)**, needed for *b* and *c* to be integers, leads to:
2n$^2$/p$^{2k}$ ∈ ℕ ⇔ 2n$^2$/p$^{2k}$ = s, s ∈ ℕ ⇔ 2n$^2$ = s·p$^{2k}$. Since p is a prime > 2: ¬(2 | p), and so (2 | s).
So s = 2t, t ∈ ℕ, and n$^2$/p$^{2k}$ = t ⇔ (n/p$^k$)$^2$ ∈ ℕ ⇔ n/p$^k$ ∈ ℕ, because n/p$^k$ is rational, so its square can only be an integer if it is an integer itself. So n/p$^k$ = r, r ∈ ℕ ⇔ **n = p$^k$·r, r ∈ ℕ, r ≥ 1.**   [63]
With **d = p$^{2k}$** and **n = p$^k$·r**, the requirement **2n$^2$ > d$^2$**, needed for *a* to be smaller than *b*, leads to:
2·(p$^k$·r)$^2$ > d$^2$ ⇔ 2·p$^{2k}$·r$^2$ > d$^2$ ⇔ 2·d r$^2$ > d$^2$ ⇔ r$^2$ > d/2 ⇔ r > $\sqrt{d/2}$
**r ≥ ceiling($\sqrt{d/2}$) = r$_{min,p,k}$**  [64]

Table T3 lists values of r and *n* for a few values of p and k. It can be easily extended.

| p | k | d = p$^{2k}$ | r$_{min,p,k}$ | r | n = p$^k$·r |
|---|---|---|---|---|---|
| 3 | 0 | 1 | 1 | 1, 2, 3, 4,.... | 1, 2, 3, 4, … |
| 3 | 1 | 9 | 3 | 3, 4, 5, 6, … | 9, 12, 15, 18, … |
| 3 | 2 | 81 | 7 | 7, 8, 9, 10, … | 63, 72, 81. 90, … |
| 3 | 3 | 729 | 20 | 20, 21, 22, 23, … | 540, 567, 594, 621, … |
| … | … | … | … | | … |
| 5 | 0 | 1 | 1 | 1, 2, 3, 4, … | 1, 2, 3, 4, … |
| 5 | 1 | 25 | 4 | 4, 5, 6, 7, … | 20, 25, 30, 35, … |
| 5 | 2 | 625 | 18 | 18, 19, 20, 21, … | 450, 475, 500, 525, … |
| … | … | … | … | | … |
| 7 | 1 | 49 | 5 | 5, 6, 7. 8. … | 35, 42, 49, 56, … |
| 7 | 2 | 2401 | 35 | 35, 36, 37, 38, … | 1715, 1764, 1813, 1862, … |
| … | … | … | … | | … |

*Table T3: values of n to form a DPT for given d = p$^{2k}$*

Note, that if k = 0, then it doesn't matter what p is, since p$^0$ = 1 for any p.

$d = p^{2k}$,  $n = p^k \cdot r$,  $2n^2/d = 2(p^k \cdot r)^2/p^{2k} = 2r^2$,  $r = r_{min,p,k}$,  $r_{min,p,k} + 1$,  $r_{min,p,k} + 2$, ...  [65r]

With [65r], the DPTs ($a$, $b$, $c$) are given by:

    $a = 2n + d$          = **$2 \cdot p^k \cdot r + p^{2k}$**     [65a]

    $b = 2n + r^2$         = **$2 \cdot p^k \cdot r + 2 \cdot r^2$**     [65b]

    $c = = a + r^2 = b + d$   = **$2 \cdot p^k \cdot r + p^{2k} + 2 \cdot r^2$**     [65c]

### ii    GCD($a$, $b$, $c$) > 1 for DPTs with $d = p^{2k}$ and (p | r)

For a DPT ($a$, $b$, $c$) given by [65], **GCD($a$, $b$, $c$) > 1 if k > 0 and p divides r**.

    (p | r) $\Rightarrow$ (p | $a$), because all the terms in [65a] are divisible by p.

    (p | r) $\Rightarrow$ (p | $b$), because all the terms in [65b] are divisible by p.

    (p | r) $\Rightarrow$ (p | $c$), because all the terms in [65c] are divisible by p.

So if k > 0 and (p | r), GCD($a$, $b$, $c$) $\geq$ p and DPT ($a$, $b$, $c$) is reducible.

### iii    IDPTs with $d = p^{2k}$, k $\geq$ 0

This section removes the DPTs from section 8ai, which result from a value of r that is divisible by p. It will be proved in section 8aiv, that the remaining DPTs are indeed IDPTs.

As is shown in section 8aii: (k > 0) $\wedge$ (p | r) $\Rightarrow$ GCD($a$, $b$, $c$) > 1.

Therefore a necessary condition for DPT ($a$, $b$, $c$) to be irreducible, is: (k > 0) $\Rightarrow$ ¬(**p | r**)     [66]

If k = 0, [65] yields the first family of IDPTs, for any p > 2 (see also section 8b below). It was proved in chapter 3 that all those DPTs are irreducible. Therefore [66] applies only to cases with k > 0.

Incorporating [66] into table T3 yields table T4, listing values of r and $n$ for a few values of p and k. It can be easily extended.

| p | k | $d = p^{2k}$ | $r_{min,p,k}$ | r | $n = p^k \cdot r$ |
|---|---|---|---|---|---|
| 3 | 0 | 1 | 1 | 1, 2, 3, 4,.... | 1, 2, 3, 4, ... |
| 3 | 1 | 9 | 3 | 4, 5, 7, 8, ... | 12, 15, 21, 24, ... |
| 3 | 2 | 81 | 7 | 7, 8, 10, 11, ... | 63, 72, 90. 99, ... |
| 3 | 3 | 729 | 20 | 20, 22, 23, 25, ... | 540, 594, 621, 675, ... |
| ... | ... | ... | ... | | ... |
| 5 | 0 | 1 | 1 | 1, 2, 3, 4, ... | 1, 2, 3, 4, ... |
| 5 | 1 | 25 | 4 | 4, 6, 7, 8, ... | 20, 30, 35, 40, ... |
| 5 | 2 | 625 | 18 | 18, 19, 21, 22, ... | 450, 475, 525, 550, ... |
| ... | ... | ... | ... | | ... |
| 7 | 1 | 49 | 5 | 5, 6, 8. 9. ... | 35, 42, 56, 63, ... |
| 7 | 2 | 2401 | 35 | 36, 37, 38, 39, ... | 1764, 1813, 1862, 1911, ... |
| ... | ... | ... | ... | | ... |

*Table T4: values of n to form an IDPT for given $d = p^{2k}$*

$d = p^{2k}$,  $n = p^k \cdot r$,  $2n^2/d = 2r^2$,  $r \geq r_{min,p,k}$,  and **if (k > 0) then ¬(p | r)**     [67n]

With [67n], the IDPTs ($a$, $b$, $c$) are given by:

    $a = 2n + d$          = **$2 \cdot p^k \cdot r + p^{2k}$**     [67a]

    $b = 2n + r^2$         = **$2 \cdot p^k \cdot r + 2 \cdot r^2$**     [67b]

    $c = = a + r^2 = b + d$   = **$2 \cdot p^k \cdot r + p^{2k} + 2 \cdot r^2$**     [67c]

### iv  Proof that the DPTs from [67] are irreducible

It will be shown, that no prime number q can be a common divisor of *a*, *b*, and *c* given by [67]. Then GCD(*a*, *b*, *c*) = 1, and (*a*, *b*, *c*) is an IDPT.

Each step below uses the obvious property: ¬(g | x) ⇒ g ∉ {CD(*x, y*)} for any numbers x and y.
See [ND2] in appendix A for the notation used.

**Step 1: Suppose q = 2. Then q ∉ {CD(*a, b, c*)}**
    [67a]: $a = 2 \cdot p^k \cdot r + p^{2k}$
    q = 2 ⇒ q is a divisor of the first term in [67a]: $2 \cdot p^k \cdot r$.
    q = 2 ⇒ q is not a divisor of the second term in [67b]: $p^{2k}$, because prime p ≠ 2.
    Using [PD4a] from appendix A: (2 | $2 \cdot p^k \cdot r$) ∧ ¬(2 | $p^{2k}$) ⇒ ¬(2 | *a*)
    So 2 ∉ {CD(*a, b, c*)}

**Step 2: Suppose q = p. Then q ∉ {CD(*a, b, c*)}**
    [67b]: $b = 2 \cdot p^k \cdot r + 2 \cdot r^2$
    q = p ⇒ q is a divisor of the first term in [75b]: $2 \cdot p^k \cdot r$.
    q = p ⇒ q is not a divisor of the second term in [75b]: $2 \cdot r^2$,
    because ¬(p | r), and using [PD1b] from appendix A: ¬(p | r) ⇒ ¬(p | $r^2$).
    Using [PD4a] from appendix A: (p | $2 \cdot p^k \cdot r$.) ∧ ¬(p | $2 \cdot r^2$) ⇒ ¬(p | *b*).
    So p ∉ {CD(*a, b, c*)}

**Step 3: Suppose q is a prime, q > 2, q ≠ p and (q | r). Then q ∉ {CD(*a, b, c*)}**
    [67a]: $a = 2 \cdot p^k \cdot r + p^{2k}$.
    (q | r) ⇒ q is a divisor of the first term in [67a]: $2 \cdot p^k \cdot r$.
    q ≠ p ⇒ q is not a divisor of the second term in [67a]: $p^{2k}$.
    Using [PD4a] from appendix A: (q | $2 \cdot p^k \cdot r$) ∧ ¬(q | $p^{2k}$) ⇒ ¬(q | *a*)
    So q ∉ {CD(*a, b, c*)}

**Step 4: Suppose q is a prime, q > 2, q ≠ p and ¬(q | r). Then q ∉ {CD(*a, b, c*)}**
    **Step 4a: Suppose q is a prime, q > 2, q ≠ p, ¬(q | r) and also ¬(q | *a*).**
        Then q ∉ {CD(*a, b, c*)}.
    **Step 4b: Suppose q is a prime, q > 2, q ≠ p, ¬(q | r) and also (q | *a*).**
        [67c]: $c = = a + 2 \cdot r^2$
        (q | *a*) ⇒ q is a divisor of the first term in [67c]: *a*.
        ¬(q | r): Using [PD1b] from appendix A: ¬(q | r) ⇒ ¬(q | $r^2$).
        So q is not a divisor of the second term in [67c]: $2 \cdot r^2$.
        Using [PD4a] from appendix A: (q | *a*) ∧ ¬(q | $2 \cdot r^2$) ⇒ ¬(q | *c*)
        So q ∉ {CD(*a, b, c*)}

Steps 1, 2, 3, 4a and 4b cover all possible prime divisors,
therefore ∀$_q$: q ∉ {CD(*a, b, c*)} ⇔ GCD(*a, b, c*) = 1
This completes the proof that all DPTs given by [67] are irreducible.

## b     Extended first family of IDPTs

Equations [67] reduce to [11] for any prime p, if k = 0:

$d = p^{2k} \Rightarrow d = 1$, $n = p^k \cdot r \Rightarrow r = n$, and $2n^2 > d^2 \Rightarrow 2n^2 > 1 \Leftrightarrow n > 0$. [67a] – [67c] reduce to:

| | | | |
|---|---|---|---|
| $a = 2 \cdot p^k \cdot r + p^{2k}$ | $= 2n + 1$ | | [67a] → [11a] |
| $b = 2 \cdot p^k \cdot r + 2 \cdot r^2$ | $= 2n + 2n^2$ | | [67b] → [11b] |
| $c = a + 2 \cdot r^2$ | $= 2n + 1 + 2n^2$ | | [67c] → [11c] |

which define the first family of IDPTs (see chapter 3) for any $n > 0$. Therefore equations [67] are a generalization of [11] and form the **extended first family of IDPTs**.

The extended first family of IDPTs is given by theorem 3 below.
It concerns all IDPTs in which $d$ is an even power of a single prime $p > 2$.

**Theorem 3**

Let $d$, $r$ and $n$ satisfy criterion C3, consisting of C3d, C3r and C3n:

| | |
|---|---|
| [C3d]: $d = p^{2k}$, p is a prime number, $p > 2$, $k \in \mathbb{N}$, $k \geq 0$. | [C3d] |
| [C3r]: $r_{min,p,k} = \text{ceiling}(\sqrt{d/2})$.    $r \in \mathbb{N}$: $(r \geq r_{min,p,k}) \wedge \{(k > 0) \Rightarrow \neg(p \mid r)\}$ | [C3r] |
| [C3n]: $n = p^k \cdot r$ | [C3n] |

Let triangles ($a$, $b$, $c$) be given by:

| | | |
|---|---|---|
| $a = 2n + d$ | $= 2 \cdot p^k \cdot r + p^{2k}$ | [67a] |
| $b = 2n + 2n^2/d$ | $= 2 \cdot p^k \cdot r + 2 \cdot r^2$ | [67b] |
| $c = a + 2n^2/d = b + d$ | $= 2 \cdot p^k \cdot r + 2 \cdot r^2 + p^{2k}$ | [67c] |

Then all triangles ($a$, $b$, $c$) are IDPTs, and no other IDPTs with $d$ satisfying C3d exist.
Proof: section 8a proves this theorem.

All these IDPTs form the **extended first family of IDPTs.**
The first family of IDPTs, discussed in chapter 3, is the case in theorem 3 with k = 0.

### c    Examples of IDPTs in the extended first family

**p = 3, k = 0**
*d* = 1, *n* = **1, 2, 3, 4, 5, 6, 7, …**
IDPTs: (3, 4, 5), (5, 12, 13), (7, 24, 25), (9, 40, 41), (11, 60, 61), (13, 84, 85), (15, 112, 113), …

**p = 3, k = 1**
*d* =9, *n* = **12, 15, 21, 24, 30, …**
IDPTs: (33, 56, 65), (39, 80, 89), (51, 140, 149), (57, 176, 185), (69, 260, 269),…

**p = 3, k = 2**
*d* =81, *n* = **63, 72, 90, 99, …**
IDPTs: (207, 224, 305), (225, 272, 353), (261, 380, 461), (279, 440, 521), …

**p = 3, k = 3**
*d* =81, *n* = **540, 594, 621, 675, …**
IDPTs: (1809, 1880, 2609), (1917, 2156, 2885), (1971, 2300, 3029), (2079, 2600, 3329), …

**p = 5, k = 0**
Identical to the first example p = 3, k = 0. For p>2, all cases with k=0 are the same.

**p = 5, k = 1**
*d* =25, *n* = **20, 30, 35, 40, …**
IDPTs: (65, 72, 97), (85, 132, 157), (95, 168, 193), (105, 208, 233), …

**p = 5, k = 2**
*d* =625, *n* = **450, 475, 525, 550, …**
IDPTs: (1525, 1548, 2173), (1575, 1672, 2297), (1675, 1932, 2557), (1725, 2068, 2693), …

**p = 7, k = 1**
*d* = 49, *n*= **35, 42, 56, 63, …**
IDPTs: (119 120 169), (133, 156, 205), (161, 240, 289), (175, 288, 337), …

**p = 11, k =3**
*d* = 1771561, *n* = **1253802, 1255133, 1256464, ….**
IDPTs: (4279165, 4282332, 6053893), (4281827, 4288764, 6060325), (4284489 4295200 6066761), …

**p = 23, k = 1**
*d* = 279841, *n* = **198375, 198904, 199433**
IDPTs: (676591, 678000, 957841), (677649, 680560, 960401), (678707, 683124, 962965), ….

**p = 127, k = 2**
*d* = 260144641, *n* = **183951245, 183967374, …**
IDPTs: (628047131, 628050540, 888195181), (628079389, 628128420, 888273061), …

…

# 9 The general case: $d = f \cdot \prod_i p_i^{2k_i}$

Chapters 7 and 8 cover all cases in which *d* is a power of a single prime number. This chapter generalizes these results for any value of *d* that satisfies theorem 1 in [37]:

**(a, b, c) is an IDPT** $\Rightarrow$ $\neg(2 \mid e_2(d)) \land \forall_{p>2}: (2 \mid e_p(d))$ [37]

To accomplish the desired generalization, the results of chapters 7 and 8 are merged in section 9b, which yields a single theorem covering the case that there is only one prime factor in PF(*d*).

This result is extended to encompass all products of primes that satisfy [37] in section 9c. In this process, the requirement is revisited that the exponents in PF(*d*) should all be greater than 0 (see appendix A, part d). This is the case in chapter 7, in which $e_2(d) = (2k+1)$, $k \geq 0$.

In chapter 8 however, in which $e_p(d) = 2k$ for the single prime divisor p, the case k = 0 was also allowed, because this enabled the case *d* = 1 to be included in the general discussion.

But if *d* can be a product of any number of odd primes, each with an even power in PD(*d*), then ambiguities arise if $k_i = 0$ for several primes. Therefore, in section c the requirement that every exponent must be greater than 0 is restored, and the case *d* = 1 will be treated as a special case.

First, in section 9a, the notation used is explained.

### a Notation

A **prime number** is written as $p_i$, with an index that simply numbers the primes.

From appendix A, part d: the following ordering is used: **i < j $\Leftrightarrow$ $p_i$ < $p_j$**. [NP3]

The product of an arbitrary number of primes $p_i$ with their exponents $k_i$ is written as $\prod_i p_i^{k_i}$, in which the index i (i = 1, 2, … ) runs over all the primes in the product.

If a certain prime p occurs in a product **g** = $\prod_i p_i^{k_i}$, then **p is a prime divisor of g**,

which (see appendix A< part d) is notated as **p $\in$ {PD(g)}** [NP5]

## b  Merging the first and second extended families of IDPTs: $d = f \cdot p^{2k}$

There is a striking similarity between theorem 2 in section 7b and theorem 3 in section 8b, but there are also a few small differences. It is possible to merge the two theorems, with only one minor adjustment, so that a single theorem covers all the cases in which *d* is a power of a single prime. It is shown below how to do so by discussing each part of the conditions and calculations in turn.

### i  Conditions C2d and C3d → C4d

[C2d] reads: $d = 2 \cdot 2^{2k}$, $k \in \mathbb{N}$, $k \geq 0$.
[C3d] reads: $d = p^{2k}$, p is a prime number, $p > 2$, $k \in \mathbb{N}$, $k \geq 0$.
These can be merged by allowing p = 2, and by defining $d = f \cdot p^{2k}$, with f=2 if p=2 and f=1 if p>2.
[C4d]:   p is a prime number. $k \in \mathbb{N}$, $k \geq 0$. $(p = 2) \Rightarrow (f = 2)$ and $(p > 2) \Rightarrow (f = 1)$.
$\qquad d = f \cdot p^{2k}$ [C4d]
Note that [C4d] reduces to [C2d] if p=2, and to [C3d] if p>2.

### ii  Conditions C2r and C3r → C4r

[C2r]: $r_{min,2,k}$ = ceiling($\sqrt{d}$).   $r \in \mathbb{N}$: $(r \geq r_{min,2,k}) \wedge \neg(2 \mid r)$.
[C3r]: $r_{min,p,k}$ = ceiling($\sqrt{d/2}$).   $r \in \mathbb{N}$: $(r \geq r_{min,p,k}) \wedge (k > 0) \Rightarrow \neg(p \mid r))$.
These can be merged by allowing p = 2, and using f as in [C4d]:
[C4r]:   p is a prime number. $k \in \mathbb{N}$, $k \geq 0$. $(p = 2) \Rightarrow (f = 2)$ and $(p > 2) \Rightarrow (f = 1)$.
$\qquad \mathbf{r_{min,p,k}}$ = **ceiling($\sqrt{fd/2}$)**.   $r \in \mathbb{N}$: $(r \geq r_{min,p,k}) \wedge \{(k > 0) \Rightarrow \neg(p \mid r)\}$ [C4r ]
Note that [C4r] reduces to [C2r] if p=2, and to [C3r] if p>r.

### iii  Conditions C2n and C3n → C4n

[C2n]: $n = 2^k \cdot r$.
[C3n]: $n = p^k \cdot r$, with p being a prime number, p > 2.
These can be merged by allowing p=2.
[C4n]:   p is a prime number. $k \in \mathbb{N}$, $k \geq 0$. r is given by {C4r].
$\qquad \mathbf{\mathit{n} = p^k \cdot r}$ [C4n]
Note that [C4n] reduces to [C2n] if p=2, and to [C3n] if p>2.

### iv  Calculating *a*, *b* and *c*: [66] and [75] → [76]

From theorem 2:

| | | |
|---|---|---|
| $a = 2n + d$ | $= 2 \cdot 2^k \cdot r + 2 \cdot 2^{2k}$ | [62a] |
| $b = 2n + 2n^2/d$ | $= 2 \cdot 2^k \cdot r + r^2$ | [62b] |
| $c = a + 2n^2/d = b + d$ | $= 2 \cdot 2^k \cdot r + r^2 + 2 \cdot 2^{2k}$ | [62c] |

From theorem 3:

| | | |
|---|---|---|
| $a = 2n + d$ | $= 2 \cdot p^k \cdot r + p^{2k}$ | [67a] |
| $b = 2n + 2n^2/d$ | $= 2 \cdot p^k \cdot r + 2 \cdot r^2$ | [67b] |
| $c = a + 2n^2/d = b + d$ | $= 2 \cdot p^k \cdot r + 2 \cdot r^2 + p^{2k}$ | [67c] |

These can be merged by allowing p=2 and using f from [c4d] and [C4r].
Then $2n^2/d = (2 \cdot p^{2k} \cdot r^2) / (f \cdot p^{2k}) = 2 \cdot r^2/f$. Then [62] and [67] merge into [68]:

| | | |
|---|---|---|
| $a = 2n + d$ | $= 2 \cdot p^k \cdot r + f \cdot p^{2k}$ | [68a] |
| $b = 2n + 2n^2/d$ | $= 2 \cdot p^k \cdot r + 2 \cdot r^2/f$ | [68b] |
| $c = a + 2n^2/d = b + d$ | $= 2 \cdot p^k \cdot r + 2 \cdot r^2/f + f \cdot p^{2k}$ | [68c] |

### v  All IDPTs with *d* being a power of a single prime: theorem 4.

**Theorem 4**

Let p be a prime number. Let k ∈ ℕ, k ≥ 0. Let f be defined by: (p=2) ⇒ (f=2) and (p>2) ⇒ (f=1).
Let *d*, r and *n* satisfy criterion C4, consisting of C4d, C4r and C4n:

| | |
|---|---|
| **d = f·p$^{2k}$** | [C4d] |
| **r$_{min,p,k}$ = ceiling($\sqrt{fd/2}$)**.  r ∈ ℕ: (r ≥ r$_{min,p,k}$) ∧ {(k > 0) ⇒ ¬(p ∣ r)} | [C4r] |
| **n = p$^k$·r** | [C4n] |

Let triangles (*a*, *b*, *c*) be given by:

| | | |
|---|---|---|
| $a = 2n + d$ | $= 2 \cdot p^k \cdot r + f \cdot p^{2k}$ | [68a] |
| $b = 2n + 2n^2/d$ | $= 2 \cdot p^k \cdot r + 2 \cdot r^2/f$ | [68b] |
| $c = a + 2n^2/d = b + d$ | $= 2 \cdot p^k \cdot r + 2 \cdot r^2/f + f \cdot p^{2k}$ | [68c] |

Then all triangles (*a*, *b*, *c*) are IDPTs, and no other IDPTs with *d* satisfying C4d exist.
Proof: Theorems 2 and 3 were proved in sections 7b and 8b.
Sections 8di – 8diii prove that criterion C4 reduces to criterion C2 if p=2, and to criterion C3 if p>2.
Section 8div proves that [68] is equivalent to [62] if p=2, and to [67] if p>2.
This proves theorem 4.

Note that the only difference between the cases p=2 and p>2 lies in the value of f:
If p = 2 then f = 2, but if p > 2 then f = 1.
This arises from the factor 2 in the term $2n^2/d$,
which interacts with prime 2 but not with the other primes.

### c    d as product of primes: $d = f \cdot \prod_i p_i^{2k_i}$

Using the results of section b, the general case can now be stated in full.

In section 9ci, this is presented in the form of theorem 5, which is proved in section 9cii, quite analogous to the way theorems 2 and 3 were proved in chapters 7 and 8.

### i    Theorem 5: necessary and sufficient conditions for (*a*, *b*, *c*) to be an IDPT.

The basis to find IDPTs is [27]: A DPT is given by

    (*a*, *b*, *c*) = ($2n + d$, $2n^2/d + 2n$, $2n^2/d + 2n + d$), with $d > 0$, $2n^2 > d^2$, ($d \mid 2n^2$)      [27]

A necessary condition for a DPT given by [27] to be an IDPT is [37]:

    (a, b, c) **is an IDPT** $\Rightarrow$ ( $\neg$(2 | $e_2(d)$) $\wedge$ $\forall_{p>2}$: (2 | $e_p(d)$))      [37]

That is, in PF(*d*), $e_2(d)$ must be odd, and $e_p(d)$ for p > 2 must be even.

Because [37] is only a necessary condition, all IDPTs satisfy [37], but there might be DPTs that also satisfy [37] that yet are not IDPTs. To discard such DPTs, a sufficient condition is also needed, which states that <u>only</u> IDPTs satisfy it. This can be accomplished by limiting the possible values of *n*, such that only values of *n* (for each *d* satisfying [37]) are allowed so that (*a*, *b*, *c*) always forms an IDPT. Theorem 5 supplies these sufficient conditions. Preliminary remarks:

In PF(*d*), all exponents must be greater than 0 (see appendix A, part d). However, the case *d* = 1 can be incorporated only if all exponents in PF(*d*) are zero. This leads to the problem, that if exponents equal to 0 are allowed in a prime factorization, then there are infinitely many ways to give a prime factorization of *d*: any extra factor $p^0$ gives the same result. This problem can be solved by the following requirement:

**If at least one exponent in PF(*d*) > 0, then all exponents in PF(d) must be greater than 0.**      [69]

So the only possibility to have $k_i = 0$, is if $k_j = 0$ for all $j \neq i$ as well.

Note, that the exponent of p = 2 in PF(*d*) cannot be 0, because $e_2(d)$ = (2k+1), and 0 cannot be written as (2k+1) since k is an integer. So condition [69] implies that 2 $\notin$ {PD(*d*)}, otherwise d > 1. This agrees with chapter 8, in which the case *d* = 1 was shown to belong to the extended first family of IDPTs, in which p ≠ 2. So [69] covers the *d* = 1 case, and prevents zero exponents for all *d* > 1.

In section 9b above, the factor f is introduced to merge the cases p = 2 and p > 2. It is needed, because the factor 2 in the term $2n^2/d$ interacts with $d = p^k$ if p = 2, but not if p > 2. In the general case, this factor is therefore needed as well, now with the condition:

(2 $\in$ {PD(*d*)} $\Rightarrow$ (f = 2)) $\wedge$ (2 $\notin$ {PD(*d*)} $\Rightarrow$ (f = 1))      [70]

Condition [C5r] in theorem 5 sets a minimum value for r: $r_{min,d}$, which depends on *d* and comes from the requirement that in (*a*, *b*, *c*): a < b. Another requirement (proved in section 9ciii) is that r cannot have a prime divisor that is also a divisor of *d*. This is written as (p $\in$ {PD(*d*)}) $\Rightarrow$ $\neg$(p | r). However, for the case *d* = 1, which allows $p_i^0$ for several $p_i$ > 2, with all exponents equal to zero, this restriction $\neg$(p | r) does not hold. Therefore this exception is incorporated in the complete condition:

([$\exists$ $k_i$ > 0] $\wedge$ [(p $\in$ {PD(*d*)}) $\Rightarrow$ $\neg$(p | r)]).

Theorem 5 generates all existing IDPTs, and only all existing IDPTs.
The proof is given in section 9ciii.

**Theorem 5**

Let numbers $d$, $r$ and $n$ satisfy the following conditions:

[C5d]     $d \in \mathbb{N}$, $d > 0$, given by $d = f \cdot \prod_i p_i^{2k_i}$. In this expression for $d$:
   (1) **$p_i$ is a prime number**.
       The index i distinguishes different primes.
       The primes in the expression for $d$ are ordered: $i < j \Leftrightarrow p_i < p_j$,
   (2) **f: ($p_1 = 2$) $\Rightarrow$ (f = 2) $\wedge$ ($p_1 > 2$) $\Rightarrow$ (f = 1).**
       Note: because $i < j \Leftrightarrow p_i < p_j$, only $p_1$ can be equal to 2.
   (3) **$k_i \in \mathbb{N}$, $k_i \geq 0$. But: ($\exists\ e_i(d) > 0$) $\Rightarrow$ ($\forall_{j \neq i}\ e_j(d) > 0$).**
       Note: if $p_i > 2$, then the total exponent $e_{p_i}(d) = 2k_i$,
       and if $p_1 = 2$, then f = 2, and the total exponent $e_2(d) = 2k_1 + 1$.

[C5r]     $r \in \mathbb{N}$, $r > 0$, satisfying the following conditions:
   (1) $r_{min,d}$ = **ceiling($\sqrt{f \cdot d/2}$)**.
   (2) $r \geq r_{min,d}$
   (3) $([\exists\ k_i > 0] \wedge [(p \in \{PD(d)\}) \Rightarrow \neg(p\ |\ r)])$

[C5n]     $n \in \mathbb{N}$, $n > 0$, satisfying the following conditions:
   (1) $n = \prod_i p_i^{k_i} \cdot r$
       Note: the factor f is not present here, and the exponents have no factor 2.

Let numbers $a$, $b$ and $c$ satisfy the following conditions:

[C5a]     $a = 2n + d$           $= 2 \cdot \prod_i p_i^{k_i} \cdot r + f \cdot \prod_i p_i^{2k_i}$

[C5b]     $b = 2n + 2n^2/d$      $= 2 \cdot \prod_i p_i^{k_i} \cdot r + 2 \cdot r^2/f$

[C5c]     $c = a + 2n^2/d = b + d$    $= 2 \cdot \prod_i p_i^{k_i} \cdot r + f \cdot \prod_i p_i^{2k_i} + 2 \cdot r^2/f$

Then all triangles ($a$, $b$, $c$) are IDPTs, and no other IDPTs exist.

         **ii**     **Proof of theorem 5**

**Step 1**     **d = 1**

Suppose $d = f \cdot \prod_i p_i^0$. Then [C5d(3)] implies $p_1 > 2$, so f = 1, d = 1, $r_{min,d} = 1$, $r \geq 1$,
and the antecedent of [C5r(3)] is false, so no divisor of r is to be considered.
r = 1, 2, 3, 4, … and n = 1, 2, 3, 4, ….
($a$, $b$, $c$) = ($2n+1$, $2n + 2n^2$, $2n + 2n^2 + 1$). This is the first family of IDPTs, see chapter 3.
It was proved in chapter 3 that all these triangles are irreducible, and since $n$ can have any integer
value, no other IDPTs with $d = 1$ exist.

**Step 2**     **Necessary condition on d**

It was proved in chapter 6, that a necessary condition for ($a$, $b$, $c$) to be an IDPT is [37]:

($a$, $b$, $c$) is an IDPT with $d > 1$  $\Rightarrow$  ($\neg(2\ |\ e_2(d))\ \wedge\ \forall_{p>2}: (2\ |\ e_p(d))$)                  [37]

Condition {C5d}: $d = f \cdot \prod_i p_i^{2k_i}$ is equivalent to [37] for $d > 1$, because
   (1) If $p_1 > 2$ then f = 1 and $d = \prod_i p_i^{2k_i}$, with $\forall_i$: ($p_i > 2$ and $(2\ |\ e_{p_i}(d))$).
       So $\forall_{p>2}: (2\ |\ e_p(d))$, and there is no $e_2(d)$, so $d = f \cdot \prod_i p_i^{2k_i}$ is equivalent to [37] if $p_1 > 2$.
   (2) If $p_1 = 2$ then f = 2 and $d = 2 \cdot 2^{2k_1} \cdot \prod_{i>1} p_i^{2k_i} = 2^{(2k_1+1)} \cdot \prod_{i>1} p_i^{2k_i}$.
       So $\neg(2\ |\ e_2(d))$ and $\forall_{p>2}: (2\ |\ e_p(d))$, so $d = f \cdot \prod_i p_i^{2k_i}$ is equivalent to [37] if $p_1 = 2$.
   (3) (1) and (2) above cover all possible cases, so {C5d} is equivalent to [37].

**Step 3**  Condition on *r*: sufficient condition for (*a*, *b*, *c*) to be irreducible

It remains to be proven that:

a  If $d > 1$, and $\exists\, p \in \{PD(d)\}$: (p | r),
   then $p \in \{CD(a, b, c)\}$, so GCD(*a*, *b*, *c*) > 1 and (*a*, *b*, *c*) is reducible.

b  if $d > 1$, and $\forall\, p \in \{PD(d)\}$: ¬(p | r)), then (*a*, *b*, *c*) is irreducible.

**Step 3a**

Suppose $\exists\, p \in \{PD(d)\}$ such that (p | r). Note that (p | d) because $p \in \{PD(d)\}$.
Then **(p | *a*)** because $a = 2 \cdot \prod_i p_i^{k_i} \cdot r + f \cdot \prod_i p_i^{2k_i} = 2 \cdot \prod_i p_i^{k_i} \cdot r + d$.
The first term is divisible by p because (p | r) and the second term is divisible by p because (p | d).
Also **(p | *b*)** because $b = 2 \cdot \prod_i p_i^{k_i} \cdot r + 2 \cdot r^2 / f$ and (p | r), so both terms are divisible by p.
Also **(p | *c*)** because $c = b + d$, and (p | b) and (p | d), so both terms are divisible by p.
So GCD(*a*, *b*, *c*) ≥ p > 1, and (*a*, *b*, *c*) is reducible.

**Step 3b**

Suppose [(p ∈ {PD(*d*)}) ⇒ ¬(p | r)]) holds, so **∀ $p_i$ ∈ {PD(*d*)}: ¬($p_i$ | r).**

i  **Suppose q is a prime number, q = $p_i$ ∈ {PD(*d*)} for some value of i. Then q ≠ {CD(*a*, *b*, *c*)}.**
   $b = 2 \cdot \prod_i p_i^{k_i} \cdot r + 2 \cdot r^2 / f$.
   q divides the first term of *b*, because by assumption q = $p_i$ for some value of i.
   q = $p_i$ does not divide the second term of *b*, because by assumption ∀ $p_i$ ∈ {PD(*d*)}: ¬($p_i$ | r),
     and using [PD1b] from appendix A: ¬($p_i$ | r) ⇔ ¬($p_i$ | $r^2$).
     If q = 2, then f = 2, and the second term of *b* is equal to $r^2$, so ¬(q | $2 \cdot r^2 / f$).
       If q > 2, and f = 2, then the second term of *b* is equal to $r^2$, so ¬(q | $r^2$).
       If q > 2, and f = 1, then the second term of *b* is equal to $2r^2$, so ¬(q | $2 \cdot r^2$).
   So in all cases, q does not divide the second term of *b*. Using [PD4a] from appendix B:
   (q | $2 \cdot \prod_i p_i^{k_i} \cdot r$) ∧ ¬(q | $2 \cdot r^2 / f$) ⇒ ¬(q | b) ⇒ q ≠ {CD(*a*, *b*, *c*)}.

ii  **Suppose q is a prime number, q ≠ $p_i$ ∈ {PD(*d*)} for any value of i.**
    **Suppose also that (q | r). Then q ≠ {CD(*a*, *b*, *c*)}.**
    $a = = 2 \cdot \prod_i p_i^{k_i} \cdot r + f \cdot \prod_i p_i^{2k_i}$.
    q divides the first term of *a*, because by assumption (q | r).
    q does not divide the second term of *a*, because:
       If q = 2, then f = 1 since by assumption q ∉ {PD(*d*)}.
         Because q ≠ $p_i$ for any value of i, ¬(q | $\prod_i p_i^{2k_i}$).
           If q > 2, and f = 2, then the second term of *a* is equal to $2 \cdot \prod_i p_i^{2k_i}$.
           Because q ≠ 2, and q ≠ $p_i$ for any value of i, ¬(q | $2 \cdot \prod_i p_i^{2k_i}$).
           If q > 2, and f = 1, then the second term of *a* is equal to $\prod_i p_i^{2k_i}$.
           Because q ≠ $p_i$ for any value of i, ¬(q | $\prod_i p_i^{2k_i}$).
    So in all cases, q does not divide the second term of *a*. Using [PD4a] from appendix B:
    (q | $2 \cdot \prod_i p_i^{k_i} \cdot r$) ∧ ¬(q | $f \cdot \prod_i p_i^{2k_i}$) ⇒ ¬(q | b) ⇒ q ≠ {CD(*a*, *b*, *c*)}.

**iii**     **Suppose q is a prime number, q ≠ $p_i$ ∈ {PD($d$)} for any value of i.**
        **Suppose also that ¬(q | r). Then q ∉ {CD($a$, $b$, $c$)}.**

     **iiia**    **Suppose q is a prime number, q ≠ $p_i$ ∈ {PD($d$)} for any value of i.**
        **Suppose also that ¬(q | r) and ¬(q | $a$). Then q ∉ {CD($a$, $b$, $c$)}.**
¬(q | $a$)   ⇒   q ∉ {CD($a$, $b$, $c$)}.
     **iiib**    **Suppose q is a prime number, q ≠ $p_i$ ∈ {PD($d$)} for any value of i.**
        **Suppose also that ¬(q | r) and (q | $a$). Then q ∉ {CD($a$, $b$, $c$)}.**
$c = a + 2 \cdot r^2 / f$.
q divides the first term of $c$ because by assumption (q | $a$).
q does not divide the second term of $c$, because:
       q = $p_i$ does not divide the second term of $b$, because by assumption ¬(q | r),
       and using [PD1b] from appendix A: ¬(q | r) ⇔ ¬(q | $r^2$).
       If q = 2, then f = 2, and the second term of $c$ is equal to $r^2$, so ¬(q | $r^2$).
          If q > 2, and f = 2, then the second term of $b$ is equal to $r^2$, so ¬(q | $r^2$).
          If q > 2, and f = 1, then the second term of $b$ is equal to $2r^2$, so ¬(q | $2r^2$).
So in all cases, q does not divide the second term of $c$. Using [PD4a] from appendix B:
(q | $2 \cdot \prod_i p_i^{k_i} \cdot r$) ∧ ¬(q | $2 \cdot r^2 / f$)   ⇒   ¬(q | $c$)   ⇒   q ∉ {CD($a$, $b$, $c$)}.

Steps 3bi – 3biii cover all possible values of q, and so ($a$, $b$, $c$) have no common factor.
GCD($a$, $b$, $c$) = 1, therefore all IDPTs given by [C5a]-[C5c] are irreducible.

Steps 1, 2 and 3 above prove theorem 5.

Theorem 5 yields only all existing IDPTs, and is therefore a constructive completeness theorem.

**All triangles ($a$, $b$, $c$) that satisfy theorem 5 are IDPTs, and no other IDPTs exist.**
**Therefore theorem 5 generates all possible IDPTs.**

# 10  R code for the general case

The following functions in R implement theorem 5 from chapter 9.
There are functions that calculate the Greatest Common Divisor of 2 or 3 numbers,
and functions that test whether a given trio (*a*, *b*, *c*) is an DPT or an IDPT.
Because the Prime Factorization (PF) of a number is needed, there are also functions that generate primes (up to a specified number), decompose a given number into its PF, or reconstruct a number from its PF. The PF is a data frame with vectors Prime and Power, with Power[i] the exponent of Prime[i] in the PF.
Finally there is a function that generates all IDPTs for a given *d* and for a specified number of *n*'s.

## a  Small tools

```
GCD <- function(p,q)
# Returns Greatest Common Divisor of p and q.
 {r <- p%%q
  if (r==0) return(q) else return(GCD(q,r))}

GCD3 <- function(a,b,c)
# Returns Greatest Common Divisor of a, b and c.
{return(GCD(GCD(a,b),c))}

IsDPT <- function(a,b,c)
# Returns T if a, b and c form a DPT.
{if (((c^2)-(a^2)-(b^2))==0) return(T) else return(F)}

IsIDPT <- function(a,b,c)
# Returns T if a, b and c form an IDPT.
{if ((IsDPT(a,b,c) & (GCD3(a,b,c)==1))) return(T) else return(F)}
```

## b    Tools for calculating a Prime Factorization

```
GeneratePrimesUpTo <- function(n_max)
 # Generates all prime numbers <= n_max
 # Returns a data frame with columns Index and Prime.
 # To find the Prime Factorization of a number n,
 # generate primes up to sqrt(n).
{Primes <- data.frame(Index=1, Prime=2)
 Primes[2,] <- c(2,3)
 n <- 3; NextRow <- 2
 repeat
  {n <- n+1
   i <- 0
   Pr_ceiling <- floor(sqrt(n))
   Pr_max <- Primes$Prime[max(Primes$Index[Primes$Prime<=Pr_ceiling])]
   i_max <- Primes$Index[Primes$Prime==Pr_max]
   repeat
    {i <- i+1
     if ((n %% Primes$Prime[Primes$Index==i]==0)|(i>i_max)) break
    }
   if (i > i_max)
    {NextRow <- NextRow + 1
     NewRow <- c(NextRow,n)
     Primes[NextRow,] <- NewRow}
     if (n >= n_max) break
    }
   return(Primes)
}
```

```
PrimeFactorsOf <- function(d)
 # Returns the prime factorization of d
 # as a data frame with columns Prime and Power.
 # All powers are greater than 0.
 # The function needs a global data frame Primes containing all primes up to d.
 # Use GeneratePrimesUpTo(sqrt(d)) to create data frame Primes first.
{PF <- data.frame(Prime=integer(), Power=integer())
 index <- 1
 repeat {while (d%%Primes$Prime[index] != 0) index <- index+1
       NewPrime <- Primes$Prime[index]
       NewPower <- 0
       while ((d%%NewPrime ==0) & (d > 1))
          {NewPower <- NewPower+1
           d <- d/NewPrime
          }
       PF[nrow(PF)+1,] <- c(NewPrime, NewPower)
       if (d==1) return(PF)
       }
}

NumFromPrimFacts <- function(PF)
 # Returns the number of which PF is the prime factorization.
 # PF must be a data frame with vectors Prime and Power.
 # Power[i] is the power of Prime[i] in the Prime Factorization.
{return(prod(PF$Prime ^ PF$Power))}
```

### c   Tools for generating IDPTs for a given *d* and max_num

```
IDPTsBruteForce <- function(d, max_num)
{IDPTs <- data.frame(d=integer(), n=integer(),
            a=integer(), b=integer(), c=integer())
 n <- 1; i <- 1
 repeat{while (((2*(n^2))%%d) != 0) n <- n+1
     a <- (2*n) + d
     b <- (2*n) + ((2*(n^2))/d)
     c <- b + d
     if (IsIDPT(a,b,c) & (!(a>b)))
       {IDPTs[nrow(IDPTs)+1,] <- c(d,n,a,b,c)
        i <- i+1
        if (i > max_num) return(IDPTs)
       }
     n <- n+1
     if (n > 1000000) return("No n under 1000000.")
    }
}
```

```r
IDPTsFromTheorem5 <- function(d, max_num)
{IDPTs <- data.frame(d=integer(), n=integer(),
            a=integer(), b=integer(), c=integer())
 # If d=1 return the first family of IDPTs
 if (d==1)
   {n <- 0
    for (index in 1:max_num)
      {n <- n+1
       a <- 2*n + 1
       b <- 2*n + 2*(n^2)
       c <- b + 1
       IDPTs[nrow(IDPTs)+1,] <- c(d,n,a,b,c)
      }
    return(IDPTs)
   }
 # If d>1 calculate IDPTs according to theorem 5
 # PFd$k contains the k-values of all the exponents
 PFd <- PrimeFactorsOf(d)
 # Check if d satisfies [C5d] for d>1
 if (((PFd$Prime[1]==2) &
   (((PFd$Power[1]%%2)==0) | any((PFd$Power[-1]%%2)!=0)))
    |
   ((PFd$Prime[1]!=2) & any((PFd$Power%%2)!=0)))
    return(cat("There are no IDPTs for d = ", d, ".", sep=""))
 # Calculate k, f, r and n
 Kfact <- PFd$Power
 Kfact[1] <- if (PFd$Prime[1]==2) (PFd$Power[1]-1)/2 else PFd$Power[1]/2
 Kfact[-1] <- PFd$Power[-1]/2
 PFd <- cbind(PFd,Kfact)
 f <- if (d%%2==0) 2 else 1
 r <- ceiling(sqrt(f*d/2))-1
 Nfact <- prod((PFd$Prime^PFd$Kfact))
 for (index in 1:max_num)
   {r<- r+1
    while (any((r%%PFd$Prime)==0)) {r <- r+1}
    n <- Nfact*r
    a <- 2*n + d
    b <- 2*n + 2*(n^2)/d
    c <- b + d
    IDPTs[nrow(IDPTs)+1,] <- c(d,n,a,b,c)
   }
 return(IDPTs)
}
```

# Appendix A: Notations and properties of divisibility

All numbers in this appendix are integers greater than zero.

### a     Notation for divisibility, common divisors and greatest common divisor

If p·q = r, then p divides r and q divides r. This is notated as:       **p·q = r** ⇔ **(p | r)** ∧ **(q | r)**.    [ND1]
If (f | p) and (f | q), then f is a Common Divisor of p and q.
Let {CD(p, q)} be the set of all common divisors of p and q. Then:        **f ∈ {CD(p, q)}**.    [ND2]
If g is the Greatest Common Divisor of p and q, then this is notated as:        **g = GCD(p, q)**.    [ND3]

### b     Floor and ceiling of a number

Suppose r ∈ ℝ is any real number.
The **floor** of r is **the largest integer** ≤ r. Notation: **floor(r)**.    [ND4]
Examples: floor(3.7) = 3;   floor(3) = 3.
The **ceiling** of r is **the smallest integer** ≥ r. Notation: **ceiling(r)**.    [ND5]
Example: ceiling(3.7) = 4;   ceiling(3) = 3.

### c     Notation for quotient and remainder of a Euclidian division

The Euclidian division of two arbitrary integers is written as **p/q = (quotient, remainder)**    [ND6]
Example: 18/7 = (2, 4).
The quotient of a Euclidian division p/q is given by the floor function:
**quotient = floor(p/q)**.    [ND7]
Example: quotient(18/7) = floor(18/7) = 2.
The remainder of a Euclidian division p/q is given by the modulo operator:
**remainder = (p mod q)**.    [ND8]
Example: remainder(18/7) = (18 mod 7) = 4.

### d     Notation for the prime factorization of a number

Every number g ∈ ℕ, g > 1 can be written uniquely as a product of powers of primes. For example: PF(2352) = $2^4 \cdot 3^1 \cdot 7^2$. This is called the **Prime Factorization of a number g: PF(g)**.    [NP1]
Note that $p^0 = 1$ for any prime p, and multiplication by 1 leaves any number unchanged. Therefore it is conventional to require, without loss of generality, that all the exponents in the prime factorization of a number g > 1 are greater than 0:

     **PF(g) = $\prod_i p_i^{k_i}$, with $k_i > 0$** in which the index i (i = 1, 2, …) simply numbers the primes.    [NP2]
To have a strict ordering in PF(g), the following convention is also used (without loss of generality):
     In PF(g): **i < j ⇔ $p_i$ < $p_j$**.    [NP3]
An individual exponent is referred to as: $e_7(2352)$ = 2: the exponent of 7 in PF(2352) is 2.
     In general: $e_{p_i}(\prod_i p_i^{k_i})$ = **$k_i$**    [NP4]
If **p is a prime divisor of g**, then this is notated as **p ∈ {PD(g)}**    [NP5]

### e     Properties of prime factorizations

The following useful properties hold for all $e_p(n)$ in PF(n), which are easy to verify.
∀$_p$: **$e_p(n \cdot m) = e_p(n) + e_p(m)$**    [PP1]
∀$_p$: **$e_p(n^m) = m \cdot e_p(n)$**    [PP2]

### f     Properties of divisibility of squares

Suppose p is a prime number, and n is any integer $\geq 0$. Then **(p | $n^2$) $\Leftrightarrow$ (p | n)**       [PD1a]

Proof: **(g | n) $\Rightarrow$ (g | $n^2$) holds for any number g**, and so also for any prime:

     (g | n) $\Leftrightarrow$ n/g = r, r $\in \mathbb{N}$. Then $n^2$/g = n·(n/g) = n·r $\in \mathbb{N}$, so (g | $n^2$).

     However, (g | $n^2$) $\Rightarrow$ (g | n) does not hold in general. Example: (9 | 144) but $\neg$(9 | 12).

     But **if g is a prime number, then (p | $n^2$) $\Rightarrow$ (p | n)**.

         If p is a prime number, then the only exponent greater than 0 in PF(p) is $e_p(p) = 1$.

         Property [P4] from section 6a states: (m | n) $\Leftrightarrow \forall_p$: $e_p(m) \leq e_p(n)$.

         So (p | $n^2$) $\Leftrightarrow e_p(p) \leq 2e_p(n)$. Since $e_p(p) = 1$, $2e_p(n) \geq 1 \Rightarrow e_p(n) \geq ½ \Rightarrow e_p(n) \geq 1$.

         Then from [P4], with $e_p(p) = 1$ and $e_p(n) \geq 1$: $e_p(p) \leq e_p(n) \Leftrightarrow$ (p | n).

Suppose p is a prime number, and n is any integer $\geq 0$. Then $\neg$**(p | $n^2$) $\Leftrightarrow \neg$(p | n)**       [PD1b]

Proof: [PD1b] is the negation of [PD1a].

### g     Properties of divisibility of sums of numbers

In this section, g is an integer > 0, and all other numbers are integers $\geq 0$.

**(g | x) $\wedge$ (g | y) $\Rightarrow$ (g | (x+y))**       [PD2a]

Proof: (g | x) $\Leftrightarrow$ x/g = p $\in \mathbb{N}$.

     (g | y) $\Leftrightarrow$ y/g = q $\in \mathbb{N}$.

     (x+y)/g = x/g + y/g = (p+q) $\in \mathbb{N} \Rightarrow$ (g | (x+y)).

Note that the reverse does not hold in general: (3 | (5+7)) does not imply (3 | 5) or (3 | 7).

**$\neg$(g | (x+y)) $\Rightarrow \neg$(g | x) $\vee \neg$(g | y)**       [PD2b]

Proof: [PD2b] is the negation of [PD2a].

**(g | x) $\wedge$ (g | (x+y)) $\Rightarrow$ (g | y)**       [PD3a]

Proof: (g | x) $\Leftrightarrow$ x/g = p $\in \mathbb{N}$.

     (g | (x+y)) $\Leftrightarrow$ (x+y)/g = q $\in \mathbb{N}$.

     (x+y)/g = x/g + y/g = p + y/g = q $\in \mathbb{N} \Rightarrow$ y/g = r $\in \mathbb{N} \Rightarrow$ (g | y)

**$\neg$(g | y) $\Rightarrow \neg$(g | x) $\vee \neg$(g | (x+y))**       [PD3b]

Proof: [PD3b] is the negation of [PD3a].

**(g | x) $\wedge \neg$(g | y) $\Rightarrow \neg$(g | (x+y))**       [PD4a]

Proof: $\neg$(g | y) $\Rightarrow \neg$(g | x) $\vee \neg$(g | (x+y))    [PD3b]

     Applying the general logical equivalence {(A $\Rightarrow$ B $\vee$ C) $\wedge$ (A $\wedge \neg$B)} $\Rightarrow$ C to [PD3b]:

     $\neg$(g | y) $\wedge \neg\neg$(g | x) $\Rightarrow \neg$(g | (x+y)) which upon rearranging yields [PD4a].

**(g | (x+y)) $\Rightarrow$ (g | x) $\vee \neg$(g | y)**       [PD4b]

Proof: [PD4b] is the negation of [PD4a].

**(g | (x+y)) $\wedge \neg$(g | x) $\Rightarrow \neg$(g | y)**       [PD5a]

Proof: (g | (x+y)) $\Rightarrow$ (g | x) $\vee \neg$(g | y)    [PD4b]

     Applying the general logical equivalence {(A $\Rightarrow$ B $\vee$ C) $\wedge$ (A $\wedge \neg$B)} $\Rightarrow$ C to [PD4b]:

     (g | (x+y)) $\wedge \neg$(g | x) $\Rightarrow \neg$(g | y).

**(g | x) $\Rightarrow$ (g | y) $\vee \neg$(g | (x+y))**       [DP5b]

Proof: The negation of [PD4a] is: (g | y) $\Rightarrow$ (g | x) $\vee \neg$(g | (x+y)).

     Swapping x and y yields [DP5b].